\documentclass[11pt]{article}
\usepackage{amsfonts}
\usepackage{amsmath}
\usepackage{amssymb}
\def\w.dim{{\rm w.dim}}

\newtheorem{thm}{Theorem}[section]

\newtheorem{lem}[thm]{Lemma}

\newtheorem{rem}[thm]{Remark}
\newtheorem{exmp}[thm]{Example}

\begin{document}
\baselineskip=15pt
\begin{center}
{\small Lecture Notes in Pure and Applied in Mathematics - Dekker
231 (2002) 301-311}
\end{center}
\vspace{1cm}

\noindent{\Large\bf Trivial Extensions of Local Rings and a \\
Conjecture of Costa} \vspace{1cm}

\noindent{\large S. Kabbaj$^{1,}$\footnote{This project has been
funded by King Fahd University of Petroleum \& minerals under
Project \# MS/TRIVIAL/224} and N. Mahdou$^2$}
\bigskip

\noindent{\small $^1$Department of Mathematics, KFUPM, P.O. Box
5046, Dhahran 31261, Saudi Arabia

\noindent$^2$Department of Mathematics, FST Fez-Sa\" iss, B.P.
2202, Univ. of Fez, Fez, Morocco} \vspace{4cm}

\noindent{\bf Abstract.} This paper partly settles a conjecture of
Costa on $(n, d)$-rings, i.e., rings in which $n$-presented
modules have projective dimension at most $d$. For this purpose, a
theorem studies the transfer of the $(n, d)$-property to trivial
extensions of local rings by their residue fields. It concludes
with a brief discussion -backed by original examples- of the
scopes and limits of our results.

\newpage
{\bf\Large 0. Introduction}\bigskip

All rings considered in this paper are commutative with identity
elements and all modules are unital. For a nonnegative integer
$n$, an $R$-module $E$ is $n$-presented if there is an exact
sequence $F_n\rightarrow F_{n-1} \rightarrow ...\rightarrow F_0
\rightarrow E\rightarrow 0$ in which each $F_i$ is a finitely
generated free $R$-module (In [1], such $E$ is said to have an
$n$-presentation). In particular, ``$0$-presented" means finitely
generated and ``$1$-presented" means finitely presented. Also,
$pd_{R} E$ will denote the projective dimension of $E$ as an
$R$-module.

In 1994, Costa [2] introduced a doubly filtered set of classes of rings throwing a brighter light on
the structures of non-Noetherian rings. Namely, for nonnegative
integers $n$ and $d$, a ring $R$ is an $(n, d)$-ring if every $n$-presented $R$-module
has projective dimension at most $d$. The Noetherianness deflates the $(n, d)$-property
to the notion of regular ring. However, outside Noetherian settings,
the richness of this classification resides in its ability to unify classic concepts such as
von Neumann regular, hereditary/Dedekind, and semi-hereditary/Pr\"ufer rings.
Costa was motivated by the sake of a deeper understanding of what makes a Pr\"ufer
domain Pr\"ufer. In this context, he asked ``what happens if we assume only that every
finitely presented (instead of generated) sub-module of a finitely generated free module
is projective?" It turned out that a non-Pr\"ufer domain having this
property exists, i.e., (In the $(n, d)$-jargon) a $(2, 1)$-domain which is not
a $(1, 1)$-domain. This gave rise to the theory of $(n, d)$-rings. Throughout,
we assume familiarity with $n$-presentation, coherence, and basics of the $(n, d)$-theory
as in [1, 2, 3, 6, 7, 8, 10].

Costa's paper [2] concludes with a number of open problems and conjectures, including
the existence of $(n, d)$-rings, specifically whether: ``{\sl There are examples of
$(n, d)$-rings which are neither $(n, d-1)$-rings nor $(n-1, d)$-rings, for all nonnegative
integers $n$ and $d$}". Some limitations are immediate; for instance, there are no (n, 0)-domains
which are not fields. Also, for $d=0$ or $n=0$ the conjecture reduces
to ``$(n, 0)$-ring not $(n-1, 0)$-ring" or ``$(0, d)$-ring not $(0, d-1)$-ring", respectively.

Let's summarize the current situation. So far, solely the cases $n \leq 2$
and $d$ arbitrary were gradually solved
in [2], [3], and [14]. These partial results were obtained using various pullbacks.
For obvious reasons, these were no longer useful for the
specific case $d=0$. Therefore, in [14], the author appealed to trivial extensions
of fields $k$ by infinite-dimensional $k$-vector spaces, and hence constructed
a $(2, 0)$-ring (also called $2$-von Neumann regular ring) which is not a $(1, 0)$-ring
(i.e., not von Neumann regular). This encouraged further work for other trivial extension contexts.

Let $A$ be a ring and $E$ an $A$-module. The trivial ring extension of $A$ by $E$ is the ring
$R = A \propto~E$ whose underlying group is $A \times E$ with multiplication given by
$(a, e)(a', e') = (aa', ae'+a'e)$. An ideal $J$ of $R$ has the form $J = I \propto E'$,
where $I$ is an ideal of $A$ and $E'$ is an $A$-submodule of $E$ such that $IE \subseteq E'$.
Considerable work, part of it summarized in Glaz's book [10] and Huckaba's book [11], has been
concerned with trivial ring extensions. These have proven to be useful in solving many
open problems and conjectures for various contexts in (commutative and non-commutative)
ring theory. See for instance [4, 5, 9, 12, 13, 15, 16, 17].

Costa's conjecture is still elusively outstanding. A complete
solution (i.e., for all nonnegative integers $n$ and $d$) would
very likely appeal to new techniques and constructions. Our aim in
this paper is much more modest. We shall resolve the case ``$n=3$
and $d$ arbitrary". For this purpose, Section 1 investigates the
transfer of the $(n, d)$-property to trivial extensions of local
(not necessarily Noetherian) rings by their residue fields. A
surprising result establishes such a transfer and hence enables us
to construct a class of $(3, d)$-rings which are neither $(3,
d-1)$-rings nor $(2, d)$-rings, for $d$ arbitrary. Section 2 is
merely an attempt to show that Theorem 1.1 and hence Example 1.4
are the best results one can get out of trivial extensions of
local rings by their residue fields.

\section{Result and Example}

\indent This section develops a result on the transfer of the $(n,
d)$-property for a particular context of trivial ring extensions,
namely, those issued from local (not necessarily Noetherian) rings
by their residue fields. This will enable us to construct a class
of $(3, d)$-rings which are neither $(3, d-1)$-rings nor $(2,
d)$-rings, for $d$ arbitrary.
\bigskip

The next theorem not only serves as a prelude to the construction
of examples, but also contributes to the study of the homological
algebra of trivial ring extensions.

\begin{thm}
Let $(A,M)$   be a local ring and let
$R =A \propto A/M$  be the trivial ring extension of $A$  by
$A/M.$  Then \\
1) $R$ is a $(3,0)$-ring provided $M$ is not finitely generated.\\
2) $R$ is not a $(2,d)$-ring, for each integer $d \geq 0$, provided
  $M$ contains a regular element.
\end{thm}

The proof of this theorem requires the next preliminary.

\begin{lem}
Let $A$ be a ring, $I$ a proper ideal of $A$, and $R$ the trivial ring
extension  of $A$ by $A/I$. Then $pd_{R} (I \propto A/I)$ and hence
$pd_{R} (0 \propto A/I)$ are infinite.
\end{lem}

\noindent{\it Proof.} Consider the exact sequence of $R$-modules
$$0 \rightarrow  I \propto A/I \rightarrow  R \rightarrow  R/(I \propto A/I) \rightarrow 0$$
We claim that $R/(I \propto A/I)$ is not projective. Deny. Then the sequence splits. Hence,
$I \propto A/I$ is generated by an idempotent element
$(a,e) =(a,e)(a,e) =(a^{2}, {0})$. So $I \propto A/I =R(a, {0})=Aa \propto 0$,
the desired  contradiction (since $A/I \not= 0$). It follows from the above sequence
that $$pd_{R}(R/(I \propto A/I)) = 1 + pd_{R}(I \propto A/I). \eqno (1)$$
Let $(x_{i})_{i \in \Delta}$ be a set of generators of $I$
and let $R^{(\Delta)}$ be a free $R$-module. Consider the exact sequence of
$R$-modules
$$0 \rightarrow  Ker(u) \rightarrow  R^{(\Delta)} \oplus R \buildrel u
\over\rightarrow  I \propto A/I \rightarrow 0$$
\noindent where $$u((a_{i},e_{i})_{i \in \Delta}, (a_{0},e_{0})) =
\displaystyle\sum _{i \in \Delta} (a_{i},e_{i})(x_{i},  {0}) +
(a_{0},e_{0})(0,   {1}) = ( \displaystyle\sum _{i \in \Delta} a_{i}x_{i},
{a_{0}})$$ since $x_{i} \in I$ for each $i \in \Delta$. Hence, $$Ker(u) =(U \propto
(A/I)^{(\Delta)}) \oplus (I \propto A/I)$$
\noindent where $U = \{(a_{i})_{i \in \Delta} \in
A^{(\Delta)} / \displaystyle\sum _{i \in \Delta} a_{i}x_{i}=0\}$. Therefore, we have the
isomorphism of $R$-modules $I
\propto A/I \cong (R^{(\Delta)}/(U \propto (A/I)^{(\Delta)})) \oplus (R/(I \propto
A/I))$. It follows that
$$pd_{R}(R/(I \propto A/I)) \leq pd_{R}(I \propto A/I). \eqno (2)$$
Clearly, (1) and (2) force $pd_{R}(I \propto A/I)$ to be infinite.

Now the exact sequence of $R$-modules
$$0 \rightarrow  I \propto A/I \rightarrow  R  \buildrel v \over\rightarrow  0
\propto A/I \rightarrow 0,$$
where $v(a,e) =(a,e)(0,  {1}) =(0,  {a})$, easily yields $pd_{R}(0 \propto A/I) = \infty $,
completing the proof of Lemma 1.2.2. $\Box$
\bigskip

\noindent {\it Proof of Theorem 1.2.1.} 1) Suppose $M$ is not finitely generated.
Let $H_0 (\not= 0)$ be a $3$-presented $R$-module and let $(z_{i})_{i=1,\ldots ,n}$
be a minimal set of generators of $H_0$ (for some positive integer $n$).
Consider the exact sequence of $R$-modules
$$0 \rightarrow  H_1 :=Ker(u_{0}) \rightarrow  R^n \buildrel u_0
\over\rightarrow  H_0  \rightarrow 0$$
where $u_{0}((r_{i})_{i=1, \ldots ,n}) =\displaystyle\sum
_{i=1}^{n}r_{i}z_{i}$. Throughout this proof we identify $R^n$
with $A^{n} \propto (A/M)^{n}$. Our aim is to prove that $H_1 =0$.
Deny. By the above exact sequence, $H_{1}$ is a $2$-presented $R$-module.
Let $(x_{i},y_{i})_{i=1,\ldots ,m}$ be a minimal set of generators of $H_{1}$
(for some positive integer~$m$). The minimality of $(z_{i})_{i=1,\ldots ,n}$
implies that $H_{1} \subseteq M^{n} \propto (A/M)^{n}$, whence $x_{i} \in M^n$
(and $y_{i} \in (A/M)^n$) for $i=1,\ldots ,m$. Consider the exact sequence of $R$-modules
$$0 \rightarrow  H_2 :=Ker(u_{1}) \rightarrow  R^m \buildrel u_1
\over\rightarrow  H_1  \rightarrow 0$$
where $u_{1}((a_{i},e_{i})_{i})
=\displaystyle\sum _{i=1}^{m}(a_{i},e_{i})(x_{i},y_{i})
=\displaystyle\sum _{i=1}^{m}(a_{i}x_{i},  {a_{i}}y_{i})$, since $x_{i} \in M^n$ for each $i$.
Then, $H_{2} =U \propto (A/M)^m$, where $U =\{(a_{i})_{i=1, \ldots ,m} \in
A^m /  \displaystyle\sum _{i=1}^{m}a_{i}x_{i} =0$ and $\displaystyle\sum
_{i=1}^{m}  {a_{i}}y_{i} =0\}$. By the above exact sequence, $H_2$ is a finitely presented
(hence generated) $R$-module, so that (via [11, Theorem 25.1]) $U$ is a finitely generated
$A$-module. Further, the minimality of $(x_{i},y_{i})_{i=1,\ldots ,m}$ yields $ U \subseteq M^m$.
Let $(t_{i})_{i=1,\ldots ,p}$ be a set of generators of $U$ and let
$(f_{i})_{i=p+1, \ldots ,p+m}$ be a basis of the $(A/M)$-vector space $(A/M)^m$. Consider
the exact sequence of $R$-modules
$$0 \rightarrow  H_3 :=Ker(u_{2}) \rightarrow  R^{p+m} \buildrel u_2
\over\rightarrow  H_2  \rightarrow 0$$
where $$u_{2}((a_{i},e_{i})_{i})
=\displaystyle\sum _{i=1}^{p}(a_{i},e_{i})(t_{i},0) + \displaystyle\sum _{i=p+1}^{p+m}(a_{i},e_{i})(0,f_{i})
=(\displaystyle\sum _{i=1}^{p}a_{i}t_{i},\displaystyle\sum _{i=p+1}^{p+m}a_{i}f_{i}),$$
since $t_i \in M^m$ for each $i=1, \ldots ,p$ and $(f_{i})_{i}$ is a basis of the $(A/M)$-vector
space $(A/M)^m$. It follows that $H_{3} \cong (V \propto (A/M)^{p}) \oplus (M^{m} \propto (A/M)^{m})$,
where $V =\{(a_{i})_{i=1, \ldots ,p} \in A^p /  \displaystyle\sum _{i=1}^{p}a_{i}t_{i} =0\}$.
By the above sequence, $H_3$ is  a finitely generated $R$-module. Hence $M \propto A/M$
is a finitely generated ideal of $R$, so $M$ is a finitely generated
ideal of $A$ by [11, Theorem 25.1], the desired contradiction.

Consequently, $H_1 =0$, forcing $H_0$ to be a free $R$-module.
Therefore, every $3$-presented $R$-module is projective
(i.e., $R$ is a $(3, 0)$-ring).
 \medskip

2) Assume that $M$ contains a regular element $m$. We must show
that $R$ is not a $(2,d)$-ring, for each integer $d \geq 0$. Let $J =R(m,0)$
and consider the exact sequence of $R$-modules
$$0 \rightarrow  Ker(v) \rightarrow  R \buildrel v \over\rightarrow  J
\rightarrow 0$$
where $v(a,e) =(a,e)(m,  {0}) =(am,  {0})$.
Clearly, $Ker(v) =0~\propto~(A/M) =R(0,  {1})$, since $m$ is a regular
element. Therefore, $Ker(v)$ is a finitely generated ideal of $R$ and hence $J$
is a finitely presented ideal of $R$. On the other hand, $pd_{R}(Ker(v))
=pd_{R}(0 \propto A/M) = \infty$ by Lemma 1.2.2, so $pd_{R}(J) = \infty $.
Finally, the exact sequence
 $$0 \rightarrow  J \rightarrow  R \rightarrow  R/J \rightarrow 0$$
yields a $2$-presented R-module, namely $R/J$, with infinite projective dimension (i.e., R is not
a $(2,d)$-ring, for each $d \geq 0$),
completing the proof. $\Box$
\bigskip

We are now able to construct a class of $(3, d)$-rings which are neither $(3, d-1)$-rings
nor $(2, d)$-rings, for $d$ arbitrary. In order to do this, we first recall
from [14] an interesting result establishing the transfer of the $(n, d)$-property
to finite direct sums.

\begin{thm} {\rm ([14, Theorem 2.4])}
A finite direct sum
$\bigoplus_{1\leq i\leq n} A_{i}$ is an $(n, d)$-ring if and only if
so is each $A_{i}$. $\Box$
\end{thm}

\begin{exmp} \rm
Let $d$ be a nonnegative integer and $B$ a
Noetherian ring of global dimension $d$. Let $(A_{0},M)$ be a nondiscrete valuation
domain and $A = A_{0}~\propto~(A_{0}/M)$ the trivial ring extension of $A_{0}$ by
$A_{0}/M$. Let $R =A \times B$ be the direct product of $A$ and $B$. Then
$R$ is a $(3, d)$-ring which is neither a $(3, d-1)$-ring nor a $(2, d)$-ring, for $d$ arbitrary
(The case $d=0$ reduces to ``$(3, 0)$-ring not $(2, 0)$-ring").
\end{exmp}

\noindent{\it Proof.} By Theorem 1.2.1, $A$ is a $(3, 0)$-ring (also called $3$-Von Neumann regular ring)
which is not a $(2,d')$-ring for each nonnegative integer $d'$. Moreover, $R$ is a
$(3,d)$-ring by [14, Theorem 2.4] since both $A$ and $B$ are $(3,d)$-rings
(by gnomonic theorems of Costa [2]). Further, $R$ is not a $(2,d)$-ring by [14, Theorem 2.4]
(since $A$ is not a $(2,d)$-ring). Finally, we claim that $R$ is not a $(3,d-1)$-ring. Deny.
Then $B$ is a $(3,d-1)$-ring by [14, Theorem 2.4]. Hence, by [2, Theorem 2.4] $B$ is a $(0,d-1)$-ring
since $B$ is Noetherian (i.e., $0$-coherent). So that $gldim(B) \leq d-1$,
the desired contradiction. $\Box$

\begin{section}{Discussion}

This section consists of a brief discussion of the scopes
and limits of our findings. This merely is an attempt to show that Theorem 1.2.1 and hence Example 1.4 are
the best results one can get out of trivial extensions of local rings by their residue fields.

\begin{rem}\rm
In Theorem 1.2.1, the $(n, d)$-property holds for a trivial ring extension of a local ring
$(A,M)$ by its residue field sans any $(n, d)$-hypothesis on the basic ring $A$. This is the first surprise.
The second one resides in the narrow scope revealed by this (strong) result, namely $n=3$ and
$d=0$. Thus, the two assertions of Theorem 1.2.1, put together with Costa's gnomonic theorems,
restrict the scope of a possible example to $n=3$ and $d$ arbitrary.
\end{rem}

Furthermore, since in Theorem 1.2.1 the upshot is controlled solely by restrictions on $M$,
the next two examples clearly illustrate its failure in case one denies these restrictions,
namely, ``{\sl $M$ is not finitely generated}" and ``{\sl $M$ contains a
regular element}", respectively.

\begin{exmp}\rm
 Let $K$ be a field and let $A =K[[X]] = K + M$,
where $M =XA$. We claim that the trivial ring extension $R$ of $A$
by $A/M (=K)$ is not an $(n,d)$-ring, for any integers $n, d \geq 0$.
\end{exmp}

\noindent{\it Proof.}
Let's first show that $R$ is Noetherian. Let $J =I \propto E$ be a proper ideal
of $R$, where $I$ is a proper ideal of $A$ and $E$ is a submodule of the simple $A$-module $A/M$
(i.e., $E=0$ or $E=A/M$). Since $A$ is a Noetherian valuation ring, $I =Aa$ for some $a \in M$.
Let $f \in A$ such that $(a,\bar {f}) \in J$.
Without loss of generality, suppose $J \not= R(a,\bar {f})$. Let $(c,\bar {g}) \in J \setminus R(a,\bar {f})$,
where $c, g \in A$, and let $c = \lambda a$, for some $\lambda \in A$.
Then $(0,\bar{g} - \lambda \bar{f}) = (c,\bar {g}) - (a,\bar {f})(\lambda,\bar {0}) \in
J \setminus R(a,\bar {f})$, so that
we may assume $c =0$ and $\bar {g} \not= \bar {0}$, i.e., $g$ is invertible in $A$. It follows that
$(0,\bar {1})=(0,\bar {g})(g^{-1},\bar {0}) \in J$ (hence $E=A/M$)   and
$(a,\bar {0})=(a,\bar {f}) - (0,\bar {g})(g^{-1}f,\bar {0}) \in J$.
Consequently, $J= (a, \bar 0)R+(0, \bar 1)R$, whence $J$ is finitely generated, as desired.

Now, by Lemma 1.2.2, $pd_{R}(0 \propto A/M) = pd_{R}R(0, 1) = \infty$, hence $gldim(R) = \infty$. Then an application of [2, Theorem 1.3(ix)]
completes the proof. $\Box$

\begin{exmp}\rm
 Let $K$ be a field and $E$ be a $K$-vector
space with infinite rank. Let $A =K \propto E$ be the trivial ring extension
of $K$ by $E$. The ring $A$ is a local $(2,0)$-ring by [14, Theorem 3.4]. Clearly, its
maximal ideal $M =0 \propto E$ is not finitely generated and consists entirely
of zero-divisors since $(0,e)M =0$, for each $e \in E$. Let $R =A \propto (A/M)$
be the trivial ring extension of $A$ by $A/M (\cong K)$. Then $R$ is a
$(2,0)$-ring (and hence Theorem 1.2.1(2) fails because of the gnomonic
property).
\end{exmp}

\noindent{\it Proof.} Let $H$ be a $2$-presented $R$-module and let $(x_{1}, \ldots ,x_{n})$ be a minimal
set of generators of $H$. Our aim is to show that $H$ is a projective $R$-module.
Consider the exact sequence of $R$-modules
$$0 \rightarrow  Ker(u) \rightarrow  R^n \buildrel u \over\rightarrow  H \rightarrow 0$$
where $u((r_{i})_{i=1, \ldots ,n}) =\displaystyle\sum _{i=1}^{n}
r_{i}x_{i}$. So, $Ker(u)$ is a finitely presented
$R$-module with $Ker(u) =U \propto E'$, where $U$ is a
submodule of $A^n$ and $E'$ is a $K$-vector subspace of $K^n$. We claim
that $Ker(u) =0$. Deny. The minimality of $(x_{1}, \ldots ,x_{n})$ yields
$$Ker(u) = U \propto E' \subseteq (M \propto A/M)R^{n} = (M \propto A/M)^n$$
since $R$ is local with maximal ideal $M \propto A/M$. Let $(y_{i},f_{i})_{i=1, \ldots ,p}$
be a minimal set of generators of
$Ker(u)$, where $y_i \in M^n$ and $f_i \in K^n$.
Consider the exact sequence of $R$-modules
$$0 \rightarrow  Ker(v) \rightarrow  R^p \buildrel v \over\rightarrow
Ker(u) (=U \propto E') \rightarrow 0$$
where $v((a_{i},e_{i})_{i=1, \ldots ,p}) =\displaystyle\sum
_{i=1}^{p}(a_{i},e_{i})(y_{i},f_{i}) =(\displaystyle\sum _{i=1}^{p}
a_{i}y_{i},\displaystyle\sum _{i=1}^{p} {a_{i}}f_{i})$. Here too
the minimality of $(y_{i},f_{i})_{i=1, \ldots ,p}$
yields $Ker(v) \subseteq (M \propto A/M)^p$; whence,
$Ker(v) =V \propto (A/M)^p$, where  $V =\{(a_{i})_{i=1, \ldots ,p}
\in A^p /  \displaystyle\sum _{i=1}^{p}a_{i}y_{i} =0\} (\subseteq M^p)$.
By the above exact sequence, $Ker(v)$ is a finitely generated $R$-module, so that
$V$ is a finitely generated $A$-module [11, Theorem 25.1]. Now, by the exact sequence
$$0 \rightarrow  V \rightarrow  A^p \buildrel w \over\rightarrow U \rightarrow 0$$
where $w((a_{i})_{i=1, \ldots ,p}) =\displaystyle\sum _{i=1}^{p}
a_{i}y_{i}$, $U$ is a finitely presented $A$-module (since $U$ is generated by $(y_{i})_{i=1, \ldots ,p}$).
Further, $U$ is an $A$-submodule of $A^n$ and  $A$ is a $(2,0)$-ring,
then $U$ is projective. In addition, $A$
is local, it follows that $U$ is a finitely generated free $A$-module. On the other hand,
$U \subseteq M^n = (0 \propto E)^n$ , so
$(0,e)U =0$ for each $e \in E$, the desired contradiction (since $U$ has a basis). $\Box$
\end{section}


\end{document}